

\magnification=\magstep1


\hsize=14.2cm
\vsize=19.5cm
\hoffset=-0.4cm
\voffset=0.5cm

\input amssym.def
\input amssym.tex

\font\teneufm=eufm10
\font\seveneufm=eufm7
\font\fiveeufm=eufm5 \newfam\eufmfam
\def\eufm{\fam\eufmfam\teneufm}
\textfont\eufmfam=\teneufm
\scriptfont\eufmfam=\seveneufm 
\scriptscriptfont\eufmfam=\fiveeufm

\def\txt#1{{\textstyle{#1}}}
\def\scr#1{{\scriptstyle{#1}}}
\def\r#1{{\rm #1}}
\def\B#1{{\Bbb #1}}
\def\ca#1{{\cal #1}}

\font\title=cmbx12 at 14pt
\font\author=cmcsc10
\font\srm=cmr8
 at11pt

\def\rightheadline{\hfil{\srm A smoothed GPY sieve}\hfil\folio}
\def\leftheadline{\rm\folio\hfil{\srm Y. Motohashi \& J. Pintz}\hfil}
\def\emptyheadline{\hfil}
\headline{\rm\ifnum\pageno=1 \emptyheadline\else
\ifodd\pageno \rightheadline \else \leftheadline\fi\fi}

\def\firstpage{\hss{\vbox to 1cm{\vfil\hbox{\rm\folio}}}\hss}
\def\emptyfootline{\hfil}
\footline{\ifnum\pageno=1\firstpage\else
\emptyfootline\fi}

\centerline{\title A Smoothed GPY Sieve}
\vskip 1cm
\centerline{by}
\vskip 0.5cm
\centerline{\author Yoichi Motohashi \& J\'anos Pintz}
\vskip 1cm
\noindent
{\bf Abstract:} Combining the arguments developed in [2] and [7], we introduce a
smoothing device to the sieve procedure [3] of D.A. Goldston, 
J. Pintz, and C.Y. Y{\i}ld{\i}r{\i}m (see [4] for its simplified version). Our
assertions embodied in Lemmas 3 and 4 imply that an improvement of the prime
number theorem of E. Bombieri, J.B. Friedlander and H. Iwaniec [1] should give rise
infinitely often to bounded differences between primes. 
\par
To this end, a rework of the main part of [7] is developed in 
Sections 2--3; thus the present article is essentially self-contained, except
for the first section which is an excerpt from [4].
\footnote{}{\noindent
2000 Mathematics Subject Classification: Primary 11N05; Secondary 11P32
\smallskip\noindent\srm
The first author was supported by KAKENHI 15540047, and 
the second author by OTKA grants No.\ T38396, T43623, T49693 and
Balaton program. 
}
\vskip 1cm
\noindent
{\bf 1.}
Let $N$ be a parameter increasing monotonically to infinity. There
are four other basic parameters $H, R, k, \ell$ in our discussion; the last two
are integers. We impose the following conditions to them:
$$
H\ll\log N\ll\log R\le\log N, \leqno(1.1)
$$
and
$$
1\le\ell\le k\ll1.
\leqno(1.2)
$$
All implicit constants in the sequel
are possibly dependent on $k,\ell$ at most; and besides, the symbol $c$
stands for a positive constant with the same dependency, whose value may
differ at each occurrence. It suffices to have $(1.2)$, since our eventual 
aim is to look into the possibility to detect the bounded differences between
primes with a certain modification of the GPY sieve. We surmise that
such a modification might be obtained by introducing a smoothing device. The
present article is, however, only to indicate that the GPY sieve admits indeed a
smoothing; it is yet to be seen if this particular smoothing contributes to our
eventual aim.
\medskip
Let
$$
\ca{H}=\{h_1,h_2,\ldots,h_k\}\subseteq[-H,H]\cap\B{Z},\leqno(1.3)
$$
with $h_i\ne h_j$ for $i\ne j$. Let us put, for a prime $p$,
$$
\Omega(p)=\{\hbox{different residue classes among $-h(\bmod\,p)$,
$h\in\ca{H}$}\}
\leqno(1.4)
$$
and write $n\in\Omega(p)$ instead of $n\,(\bmod\,p)\in\Omega(p)$. We call  
$\ca{H}$  admissible if
$$
|\Omega(p)|<p\quad \hbox{for all $p$},\leqno(1.5)
$$
and assume this unless otherwise stated.
We extend $\Omega$ multiplicatively, so that 
$n\in\Omega(d)$ with square-free $d$ if and only if
$n\in\Omega(p)$ for all $p|d$, which is equivalent to
$$
(n+h_1)(n+h_2)\cdots(n+h_k)\equiv 0\,(\bmod\, d).
\leqno(1.6)
$$
We put, with $\mu$ the M\"obius function,
$$
\lambda_R(d;\ell)=\cases{\hfil 0 & if $d>R$,\cr
\displaystyle{\mu(d)\over (k+\ell)!}
\left(\log {R\over d}\right)^{k+\ell} & if $d\le R$,}\leqno(1.7)
$$
and
$$
\Lambda_R(n;\ca{H},\ell)
= \sum_{n\in\Omega(d)}\lambda_R(d;\ell).\leqno(1.8)
$$
\par
Also, let
$$
E^*(y;a,q)=\vartheta^*(y;a,q)-{y\over\varphi(q)},\quad
\vartheta^*(y;a,q)=\sum_{\scr{y<n\le
2y}\atop\scr{n\equiv a\,(\r{mod}\, q)}}\varpi(n),\leqno(1.9)
$$
where $\varphi$ is the Euler
totient function; and $\varpi(n)=\log n$ if $n$ is a prime, 
and $=0$ otherwise. In all of the existing accounts [2]--[4] 
of the GPY sieve, it is assumed that
$$
\sum_{q\le x^\theta}\max_{(a,q)=1}
\max_{y\le x}|E^*(y;a,q)|\ll {x\over(\log x)^{C_0}},\leqno(1.10)
$$
with a certain absolute constant $\theta\in(0,1)$ and an arbitrary fixed $C_0>0$;
the implied constant depending only on $C_0$.
\medskip
The following asymptotic formulas $(1.12)$ and $(1.14)$ are the 
implements with which Goldston, Pintz and Y{\i}ld{\i}r{\i}m established
$$
\liminf_{n\to\infty} {\r{p}_{n+1}-\r{p}_n\over \log \r{p}_n}=0,\leqno(1.11)
$$
where $\r{p}_n$ is the $n$th prime. 
\medskip
\noindent
{\bf Lemma 1.} {\it Provided $(1.1)$, $(1.2)$, and $R\le N^{1/2}/(\log
N)^C$ hold with a sufficiently large $C>0$ depending only on $k$ and $\ell$, 
we have
$$
\leqalignno{
&\sum_{N<n\le 2N}\Lambda_R(n;\ca{H},\ell)^2
&(1.12)\cr
&={{\eufm S}(\ca{H})\over(k+2\ell)!}
{2\ell\choose\ell}N(\log R)^{k+2\ell}
+O(N(\log N)^{k+2\ell-1}(\log\log N)^c),
}
$$
where 
$$
{\eufm S}(\ca{H})=\prod_p\left(1-{|\Omega(p)|\over p}\right)
\left(1-{1\over p}\right)^{-k}.\leqno(1.13)
$$
}
\medskip
\noindent
{\bf Lemma 2.} {\it Provided $(1.1)$, $(1.2)$, $(1.10)$, and
$R\le N^{\theta/2}/(\log N)^C$ hold with a sufficiently
large $C>0$ depending only on $k$ and $\ell$,
$$
\leqalignno{
&\sum_{N<n\le 2N}\varpi(n+h)\Lambda_R(n;\ca{H},\ell)^2&(1.14)\cr
&={{\eufm S}(\ca{H})\over(k+2\ell+1)!}
{2(\ell+1)\choose\ell+1}N(\log R)^{k+2\ell+1}
+O(N(\log N)^{k+2\ell}(\log\log N)^c),
}
$$
whenever $h\in\ca{H}$.
}
\medskip
\noindent
A short self-contained treatment of the
assertions $(1.11)$--$(1.14)$  can be found in [4].
\medskip
Note that the case $h\notin\ca{H}$ in the last lemma, which is included 
in [2]--[4], is irrelevant for our present purpose.  
In fact, a combination of $(1.10)$, $(1.12)$, and $(1.14)$ gives, for
$R=N^{\theta/2}/(\log N)^{C_0}$,
$$
\leqalignno{
&\sum_{N<n\le2N}\left\{\sum_{h\in\ca{H}}\varpi(n+h)-\log3N\right\}
\Lambda_R(n;\ca{H},\ell)^2&(1.15)\cr
=&(1+o(1)){{\eufm S}(\ca{H})\over(k+2\ell)!}
{2\ell\choose\ell}N(\log R)^{k+2\ell}(\log N)
\left({k\over k+2\ell+1}\cdot{2(2\ell+1)\over\ell+1}\cdot
{\theta\over2}-1\right).
}
$$ 
Thus, the $k$-tuple $(n+h_1,\ldots,n+h_k)$ with any fixed admissible $\ca{H}$
should contain two primes for infinitely many $n$, if the last factor in
$(1.15)$ is positive. Namely, with an appropriate choice of $k,\ell$ depending on
$\theta$ we would be able to conclude that 
$$
\liminf_{n\to\infty}\,(\r{p}_{n+1}-\r{p}_n)<\infty,\leqno(1.16)
$$
provided $\theta>{1\over2}$.
\medskip
The aim of the present work is to prove a smoothed
version of $(1.12)$ and $(1.14)$ in order to look into the possibility of
replacing $(1.10)$ with a $\theta>{1\over2}$ by any less stringent hypothesis.
\medskip
In passing, we note that the historical aspect of the Selberg sieve and  
the bilinear structure of its error term can be found in [8],
including that of smoothed sieves which came later, and
are naturally relevant to our present work.
\medskip
\noindent
{\it Convention\/}. All symbols and conditions introduced above are retained.
We assume additionally that 
$$
\hbox{$H=H(k,\ell)$ is bounded,}\leqno(1.17)
$$ 
which should not cause any
loss of generality under the present circumstance. Implicit constants may depend
on $k$ at most, but they can be regarded to be absolute once the least possible
value of $k$ is fixed. Thus the dependency on $k$ of estimations will not be
mentioned repeatedly, excepting at $(4.15)$, $(4.16)$, $(5.15)$, and $(6.1)$.
\bigskip
\noindent
{\bf 2.} We shall first rework the main part of [7] in the present and the next
sections (cf.\ [6, Sections 2.3 and 3.4]). 
\medskip  
Thus let us put
$$
R_0=\exp\left({\log R\over (\log\log R)^{1/5}}\right),\quad 
R_1=\exp\left({\log R\over (\log\log R)^{9/10}}\right).\leqno(2.1)
$$
We divide the half-line $(R_0,\infty)$ into intervals
$(R_0R_1^{j-1},R_0R_1^j]$, $j=1,2,\ldots$, denoting them by $P$, with or
without suffix. We let $|P|$ be the right end point of $P$. 
\medskip
Let
$$
R_0R_1\le z\le R.\leqno(2.2)
$$
We consider the commutative semi-group $\ca{Y}(z)$ generated  by all $P$
such that $|P|\le z$.
Let $D=P_1P_2\cdots P_r$  be an
element of $\ca{Y}(z)$. Then the notation $d\in D$   
indicates that $d$ has the prime decomposition
$d=p_1p_2\cdots p_r$ with
$p_j\in P_j$ ($1\le j\le r$). We use the convention $1\in D$ if and only if $D$ is
the empty product. Also,
$|D|$ stands for $|P_1|\cdots|P_r|$. Naturally, $|D|=1$ if $D$ is empty.
\medskip
Let $\xi$ be a real valued function over $\ca{Y}(z)$, which
satisfies the following conditions:
$$
\xi(D)=\cases{0, &if $|D|> R$, \cr 0, & if
$D$ is not square-free,\cr
\hbox{arbitrary},& otherwise,}\leqno(2.3)
$$
with an obvious abuse of terminology.
We are concerned with the quadratic form
$$
\ca{J}=\sum_{D_1,\,D_2}\xi(D_1)\xi(D_2)
\sum_{d_1\in D_1,\,d_2\in D_2}{|\Omega([d_1,d_2])|\over
[d_1,d_2]},\leqno(2.4)
$$
where $[d_1,d_2]$ is the least common multiple of $d_1,\,d_2$. 
\medskip
In the inner double sum 
of $(2.4)$, $D_1$ and $D_2$ can be supposed to be square-free, and 
by multiplicativity the sum is equal to
$$
\leqalignno{
&\prod_{P_1|D_1,\,P_2|D_2}\left(\sum_{p_1\in P_1,\,p_2\in P_2}
{|\Omega([p_1,p_2])|\over[p_1,p_2]}\right)&(2.5)\cr
=&\prod_{P_1|D_1}\left(\sum_{p_1\in P_1}{|\Omega(p_1)|\over p_1}\right)
\prod_{P_2|D_2}\left(\sum_{p_2\in P_2}{|\Omega(p_2)|\over p_2}\right)
\prod_{\scr{P|D_1}\atop\scr{P|D_2}}{\left(\displaystyle\sum_{p,p'\in P}
{|\Omega([p,p'])|\over[p,p']}\right)\over
\left(\displaystyle\sum_{p\in P}{|\Omega(p)|\over p}\right)^2},
}
$$
with primes $p_1,p_2, p, p'$. We then introduce
$$
\leqalignno{
\Delta(D)&=\prod_{P|D}\left(\sum_{p\in P}{|\Omega(p)|\over p}\right),&(2.6)\cr
\Phi(D)&={1\over\Delta(D)^2}\prod_{P|D}\left(\sum_{p\in P}{|\Omega(p)|\over p}
\left(1-{|\Omega(p)|\over p}\right)\right).&(2.7)
}
$$
Obviously $\Phi$ does not vanish; actually we have here $|\Omega(p)|=k$ but we
retain the notation because of a future purpose. We
have, for any square-free
$D$,
$$
\sum_{K|D}\Phi(K)={1\over\Delta(D)^2}
\prod_{P|D}\left(\sum_{p,p'\in P}{|\Omega([p,p'])|\over [p,p']}\right),
\leqno(2.8)
$$
which is to be compared with the last factor in $(2.5)$. 
\medskip
From $(2.5)$--$(2.8)$, we get the diagonalisation
$$
\ca{J}=\sum_{K}\Phi(K)\Xi(K)^2,\leqno(2.9)
$$
with
$$
\Xi(K)=\sum_{K|D}\Delta(D)\xi(D).\leqno(2.10)
$$
Note that $(2.3)$ induces the restriction that
$K$ be square-free and $|K|\le R$ in $(2.9)$.
Reversing $(2.10)$, we have,
with an obvious generalisation of the M\"obius function,
$$
\xi(D)={1\over\Delta(D)}\sum_K\mu(K)\Xi(KD);\leqno(2.11)
$$
and the case $D=\emptyset$ the empty product implies that
$$
\ca{J}=\sum_K\Phi(K)\left(\Xi(K)-{\xi(\emptyset)\over
G(R,z)}{\mu(K)\over\Phi(K)}\right)^2 
+{\xi(\emptyset)^2\over G(R,z)}\,,\leqno(2.12)
$$
where
$$
G(y,z)=\sum_{|K|\le y}{\mu(K)^2\over\Phi(K)}.\leqno(2.13)
$$
Note that the appearance of $z$ here indicates that $K\in\ca{Y}(z)$.
\medskip
We now set 
$$
\Xi(K)=\xi(\emptyset){\mu(K)\over G(R,z)\Phi(K)}\leqno(2.14)
$$ 
or by (2.11)
$$
\xi(D)={\xi(\emptyset)\over
G(R,z)}{\mu(D)\over\Delta(D)\Phi(D)}
\sum_{\scr{|K|\le R/|D|}\atop\scr{(K,D)=1}}{\mu(K)^2\over\Phi(K)}.\leqno(2.15)
$$
Then we have
$$
\ca{J}={\xi(\emptyset)^2\over G(R,z)}.\leqno(2.16)
$$
Hereafter we shall work with $(2.15)$, as $(2.3)$ is obviously satisfied.
It should be noted that we have now
$$
|\xi(D)|\le |\xi(\emptyset)|,\leqno(2.17)
$$
since 
$$
G(R,z)\ge\sum_{L|D}{\mu(L)^2\over\Phi(L)}
\sum_{\scr{|K|\le R/|D|}\atop\scr{(K,D)=1}}{\mu(K)^2\over\Phi(K)}\leqno(2.18)
$$
and by $(2.8)$
$$
\leqalignno{
\sum_{L|D}{\mu(L)^2\over\Phi(L)}
&={1\over\Phi(D)}\sum_{L|D}\mu^2(L)\Phi(L)&(2.19)\cr
&={1\over\Delta(D)\Phi(D)}\cdot{1\over\Delta(D)}
\prod_{P|D}\left(\sum_{p,p'\in P}{|\Omega([p,p'])\over [p,p']}\right)\cr
&\ge{1\over\Delta(D)\Phi(D)}.
}
$$
\bigskip
\noindent
{\bf 3.} In this section we shall evaluate $G(z)=G(z,z)$
asymptotically; we are still working with $\ca{Y}(z)$.
In fact, we shall treat more generally $G(z;Q)$ with $Q\in\ca{Y}(z)$,
$\log |Q|\ll\log R$, which is
the result of imposing the restriction $y=z$ and $(K,Q)=1$ to the sum $(2.13)$. 
\medskip
We define $G(y,z;Q)$ analogously, and introduce
$$
T(y,z;Q)=\int_1^y G(t,z;Q){dt\over t},\quad T_1(y,z;Q)=
\sum_{\scr{|K|\le y}\atop\scr{(K,Q)=1}}
{\mu(K)^2\over\Phi(K)}\log|K|,\leqno(3.1)
$$
so that
$$
G(y,z;Q)\log y=T(y,z;Q)+T_1(y,z;Q).\leqno(3.2)
$$
Observe that for $1\le y< R_0R_1$ 
$$
G(y,z;Q)=1,\quad T(y,z;Q)=\log y,\quad T_1(y,z;Q)=0.\leqno(3.3)
$$
Since $\log |K|=\sum_{P|K}\log|P|$ for square-free $K$, we have
$$
T_1(y,z;Q)=\sum_{\scr{|P|\le z}\atop\scr{P\nmid Q}}
{\log|P|\over\Phi(P)}G(y/|P|,z;PQ).\leqno(3.4)
$$
\medskip
On the other hand we see readily that for any $P\nmid Q$, $|P|\le z$,
$$
G(y,z;Q)=G(y,z;PQ)+{1\over\Phi(P)}G(y/|P|,z;PQ).\leqno(3.5)
$$
Let
$$
\hbox{$\Psi(P)=\left(1+\Phi(P)\right)^{-1}\,$ or 
$\,\Phi(P)\Psi(P)=1-\Psi(P)$,}\leqno(3.6)
$$
and rewrite $(3.5)$. In the result we replace $y$ by $y/|P|$, and get
$$
\leqalignno{
G(y/|P|,z;PQ)&=\Psi(P)\Phi(P)G(y/|P|,z;Q)&(3.7)\cr
&+\Psi(P)\left\{G(y/|P|,z;PQ)-G(y/|P|^2,z;PQ)\right\}.
}
$$
Inserting this into $(3.4)$, we have that
$$
\leqalignno{
T_1(y,z;Q)&=\sum_{\scr{|K|\le y}\atop\scr{(K,Q)=1}}
{\mu(K)^2\over\Phi(K)}\sum_{\scr{|P|\le
y/|K|}\atop\scr{P\nmid Q}}\Psi(P)\log|P|&(3.8)\cr  &
+\sum_{\scr{y/z^2< |K|\le
y}\atop\scr{(K,Q)=1}}{\mu(K)^2\over\Phi(K)}
\sum_{\scr{\sqrt{y/|K|}<|P|\le y/|K|}
\atop\scr{P\nmid KQ}}
{\Psi(P)\over\Phi(P)}\log|P|,
}
$$
where the additional condition $R_0R_1\le |P|\le z$ is implicit.
\medskip
Now, to evaluate the first sum over $P$ on the right side of $(3.8)$, we observe
that since by $(2.1)$ we have
$$
\Delta(P)\ll{\log R_1\over\log |P|},\leqno(3.9)
$$
it holds that
$$
\Psi(P)=\Delta(P)\left(1+O\left(\Delta(P)\right)\right).\leqno(3.10)
$$
This implies that for $\log (R_0R_1)\le \log x\ll \log R$
$$
\sum_{\scr{|P|\le x}\atop\scr{P\nmid Q}}\Psi(P)\log|P|=
k\log{x}+O\left(\log R_0\right),
\leqno(3.11)
$$
where the implied constant is independent of $Q$.
In fact, the left side is equal to
$$
\leqalignno{
&\sum_{|P|\le x}\Delta(P)\log|P|
+O\left(\sum_{P|Q}\log R_1+\sum_{|P|\le x}{(\log
R_1)^2\over\log|P|}\right)&(3.12)\cr =&\sum_{|P|\le x}\sum_{p\in
P}{|\Omega(p)|\over p}\left(\log p+\log(|P|/p)\right)\cr &\hskip 1cm+
O\left({\log |Q|\over\log R_0}\log R_1+
\sum_{\scr{j}\atop\scr{R_0\le R_0R_1^j\ll x}}
{(\log R_1)^2\over\log(R_0R_1^j)}\right)\cr =&k\log{x\over R_0}
+O\left({\log R\over\log R_0}\log R_1\right). }
$$
Also, by $(2.7)$ and $(3.9)$--$(3.10)$,
$$
\sum_{\scr{\sqrt{y/|K|}<|P|\le y/|K|}\atop\scr{P\nmid KQ}}
{\Psi(P)\over\Phi(P)}\log|P|\ll
\sum_{\sqrt{y/|K|}<|P|\le y/|K|}
{(\log R_1)^2\over\log |P|}\ll\log R_1.\leqno(3.13)
$$
We insert $(3.11)$ and $(3.13)$ into $(3.8)$, on noting the
implicit condition mentioned there. We see that
$$
T_1(y,z;Q)=k\sum_{\scr{|K|\le
y}\atop\scr{(K,Q)=1}}{\mu(K)^2\over\Phi(K)}
\log{y\over|K|}- k\sum_{\scr{|K|\le
y/z}\atop\scr{(K,Q)=1}}{\mu(K)^2\over\Phi(K)}
\log{y/z\over|K|}+U(y,z;Q),\leqno(3.14)
$$
with
$$
U(y,z;Q)\ll G(y,z;Q)\log R_0,\leqno(3.15)
$$
provided $\log (R_0R_1)\le\log y\ll \log R$ and $\log|Q|\ll\log R$; the implied
constant is independent of $Q$.
\medskip
We set $y=z$ in $(3.2)$ and $(3.14)$, and get
$$
G(z;Q)\log z = (k+1)T(z,z;Q)+U(z,z;Q).\leqno(3.16)
$$
We are then led to the assertion that uniformly in 
$Q\in\ca{Y}(z)$, $\log|Q|\ll\log R$,
$$
G(z;Q)={W(R_0)\over k!{\eufm S}(\ca{H})}
(\log z)^k\left(1+O\left({\log
R_0\over\log z}\right)\right),\leqno(3.17)
$$
where
$$
W(R_0)=\prod_{p\le R_0}\left(1-{|\Omega(p)|\over p}\right).\leqno(3.18)
$$
The deduction of $(3.17)$ from $(3.16)$ is standard; cf.\ [5, Section
2.2.2]. We should remark in this context that 
$$
{W(R_0)\over {\eufm
S}(\ca{H})}=\left(1+O\left({\log
R_1\over\log R_0}\right)\right)\lim_{s\to0^+}\zeta(s+1)^{-k}
\prod_{P\nmid Q}\left(1+{1\over|P|^s\Phi(P)}\right);\leqno(3.19)
$$
see [7, pp.\ 1060--1601] together with a minor correction. Here $\zeta$ is the
Riemann zeta-function. Note that the left side of $(3.19)$ is independent of $Q$.
\bigskip
\noindent
{\bf 4.} With this, we are now ready to  
start smoothing the assertions of
Lemmas 1 and 2. Hereafter we shall work with
$\ca{Y}(w)$ in place of $\ca{Y}(z)$, where
$$
w=R^\omega.\leqno(4.1)
$$
The constant $\omega$ is to satisfy
$$
3\log k\le k\omega\le\txt{1\over2}k ,\leqno(4.2)
$$
while $k$ is assumed to be sufficiently large, though bounded.
\medskip
We put
$$
\tilde\lambda_R(D;\ell)={{\eufm S}(\ca{H})\over\ell!W(R_0)}
{\mu(D)\over\Phi(D)\Delta(D)}\sum_{\scr{|K|\le R/|D|}\atop
\scr{(K,D)=1}}{\mu(K)^2\over\Phi(K)}\left(\log {R/|D|\over|K|}\right)^\ell,
\leqno(4.3)
$$
where $D,K\in\ca{Y}(w)$.
This is to be compared with $(2.15)$ specialised by $z=w$ and
$$
\xi(\emptyset)={\eufm S}(\ca{H}){G(R,w)\over\ell!W(R_0)}.\leqno(4.4)
$$
The side condition $(2.3)$ is obviously satisfied; also,
by $(2.17)$ and $(4.4)$,
$$
|\tilde\lambda_R(D;\ell)|\le
{\eufm S}(\ca{H}){G(R,w)\over\ell!W(R_0)}(\log R)^\ell
\ll (\log R)^{k+\ell},\leqno(4.5)
$$
where $(3.17)$ is used via $G(R,w)\le G(R)$. Our counterpart of $(1.8)$ is now
defined to be
$$
\tilde\Lambda_R(n;\ca{H},\ell)=
\sum_{D}\tilde\lambda_R(D;\ell)\sum_{\scr{d\in D}
\atop\scr{n\in\Omega(d)}}1.\leqno(4.6)
$$
\par
As to the interval $[1, R_0]$, which is excluded in the above, we appeal to
the Fundamental Lemma in the sieve method (see [5, p.\ 92]). Thus, there
exists a function $\varrho$, supported on the set of square-free integers, 
such that $\varrho(d)=0$ or $\pm1$ for any $d\ge1$, and $\varrho(d)=0$ either if
$d\ge R_0^\tau$ with $\tau$ to be fixed later or if $d$ has a prime factor greater
than or  equal to $R_0$, and that for any $n\ge1$
$$
\gamma(n;\ca{H})
=\sum_{n\in\Omega(d)}\varrho(d)\ge0\leqno(4.7)
$$
as well as
$$
\sum_{d}{\varrho(d)\over
d}|\Omega(d)|=W(R_0)\left(1+O(e^{-\tau})\right).
\leqno(4.8)
$$
We set
$$
\tau=(\log\log R)^{1/10}.\leqno(4.9)
$$
\medskip
Now our task is to evaluate asymptotically the sum
$$
\sum_{N<n\le 2N}\gamma(n;\ca{H})
\tilde\Lambda_R(n;\ca{H},\ell)^2,\leqno(4.10)
$$
which is to replace the left side of $(1.12)$.
By $(4.5)$ and $(4.8)$, this is equal to
$$
NW(R_0)\tilde{\ca{T}}\left(1+O(e^{-\tau})\right)
+O\left(R_0^\tau R^2(\log N)^c\right),\leqno(4.11)
$$
where $\tilde{\ca{T}}$ is defined analogously to $(2.4)$. We have
$$
\leqalignno{
\tilde{\ca{T}}&=\sum_{|D|\le R}\Phi(D)\left(\sum_{D|K}
\Delta(K)\tilde\lambda_R(K;\ell)\right)^2&(4.12)\cr
&=\left({{\eufm S}(\ca{H})\over\ell!W(R_0)}\right)^2
\sum_{|K|\le R}{\mu(K)^2\over\Phi(K)}\left(\log
{R\over|K|}\right)^{2\ell}; 
}
$$
the second line is due to the relation similar to that among $(2.10)$, $(2.14)$
and $(2.15)$.  
\medskip
The last sum is 
$$
\leqalignno{
&\le \sum_{\scr{|K|\le R}\atop\scr{P|K\Rightarrow |P|\le R
}}{\mu(K)^2\over\Phi(K)}\left(\log
{R\over|K|}\right)^{2\ell}&(4.13)\cr
&=\int^R_{R_0R_1} (\log R/t)^{2\ell}dG(t)+(\log R)^{2\ell}\cr
&={(2\ell)!\over (k+2\ell)!}{W(R_0)\over{\eufm S}(\ca{H})}(\log R)^{k+2\ell}
\left(1+O\left((\log\log R)^{-1/5}\right)\right).
}
$$
In the first line we have moved to the semi-group $\ca{Y}(R)$; the second line
depends on $G(t,R)=G(t)$ for $t\le R$, and the last on $(3.17)$ with $Q=\emptyset$.
On the other hand, the Buchstab identity implies that 
the sum in question is equal
to
$$
\leqalignno{
&\sum_{\scr{|K|\le R}\atop\scr{P|K\Rightarrow |P|\le R
}}{\mu(K)^2\over\Phi(K)}\left(\log
{R\over|K|}\right)^{2\ell}&(4.14)\cr
&-\sum_{w<|P|\le R}{1\over\Phi(P)}
\sum_{\scr{|K|\le R/|P|}\atop\scr{P'|K\Rightarrow |P'|< |P|
}}{\mu(K)^2\over\Phi(K)}\left(\log
{R/|P|\over|K|}\right)^{2\ell}.
}
$$
The last double sum is
$$
\leqalignno{
&\le \sum_{w<|P|\le R}{1\over\Phi(P)}
\sum_{\scr{|K|\le R/w}\atop\scr{P'|K\Rightarrow |P'|< R/w
}}{\mu(K)^2\over\Phi(K)}\left(\log
{R/w\over|K|}\right)^{2\ell}&(4.15)\cr
&\ll k|\log\omega|{(2\ell)!\over (k+2\ell)!}{W(R_0)\over{\eufm S}(\ca{H})}(\log
R/w)^{k+2\ell}\cr &\ll e^{-k\omega/3}{(2\ell)!\over (k+2\ell)!}
{W(R_0)\over{\eufm S}(\ca{H})}(\log R)^{k+2\ell}, }
$$
where $(4.2)$ has been invoked, and the implied constants are absolute.
\medskip
Hence collecting $(4.11)$--$(4.15)$ we obtain the following smoothed version
of Lemma 1:
\medskip
\noindent
{\bf Lemma 3.} {\it With $(1.17)$, $(2.1)$, $(4.1)$, $(4.2)$, $(4.3)$ $(4.6)$, 
$(4.7)$, $(4.9)$ and the same assumption as in Lemma 1,  we have,
as $N\to\infty$, 
$$
\leqalignno{
&\sum_{N<n\le 2N}\gamma(n;\ca{H})\tilde\Lambda_R(n;\ca{H},\ell)^2
&(4.16)\cr
&={{\eufm S}(\ca{H})\over(k+2\ell)!}
{2\ell\choose\ell}N(\log R)^{k+2\ell}\left(1+O(e^{-k\omega/3})\right),
}
$$
where the implied constant is absolute.
}
\bigskip
\noindent
{\bf 5.} Next, we shall consider a twist of $(4.16)$ with primes:
$$
\leqalignno{
&\sum_{N<n\le 2N}\varpi(n+h)\gamma(n;\ca{H})\tilde\Lambda_R(n;\ca{H},\ell)^2
&(5.1)\cr
=&\sum_{N<n\le 2N}\varpi(n+h)\gamma(n;\ca{H}\backslash\{h\})
\tilde\Lambda_R(n;\ca{H}\backslash\{h\},\ell)^2,
}
$$
as it is assumed that $h\in\ca{H}$, $R<N$.  Note that we are working with
$\ca{Y}(w)$. Expanding out the square, we see that this
is equal to
$$
\leqalignno{
&\sum_{D_1,D_2}
\tilde\lambda_R(D_1;\ell)\tilde\lambda_R(D_2;\ell)\sum_d\varrho(d)&(5.2)\cr
&\times\sum_{d_1\in D_1,d_2\in D_2}\,
\sum_{\scr{a\in\Omega^-(d[d_1,d_2])}\atop
\scr{(a+h,d[d_1,d_2])=1}}\vartheta^*(N;a+h,d[d_1,d_2])
+O(R_0^\tau R^2(\log N)^c),
}
$$
where $\Omega^-$ corresponds to $\ca{H}\backslash\{h\}$, and $(4.5)$ has been
applied. The condition in the inner-most sum induces the introduction of
$$
\Omega^*(p)=\Omega^-(p)\backslash\{-h\bmod p\}=\Omega(p)\backslash\{-h\bmod p\}.
\leqno(5.3)
$$
Note that $|\Omega^*(p)|=|\Omega(p)|-1$, which we may assume 
does not vanish, provided $p$ is sufficiently large.
\medskip
The sum in $(5.2)$ is equal to
$$
N\ca{T}^*\sum_d{\varrho(d)\over\varphi(d)}|\Omega^*(d)|+\ca{E},\leqno(5.4)
$$
where 
$$
\ca{T}^*=\sum_{D_1,D_2}\tilde\lambda_R(D_1;\ell)\tilde\lambda_R(D_2;\ell)
\sum_{d_1\in D_1,\, d_2\in
D_2}{|\Omega^*([d_1,d_2])|\over\varphi([d_1,d_2])}\leqno(5.5)
$$
and
$$
\ca{E}=\sum_{D_1,D_2}
\tilde\lambda_R(D_1;\ell)\tilde\lambda_R(D_2;\ell)\sum_d\varrho(d)
\sum_{d_1\in D_1,d_2\in D_2}\,
\sum_{a\in\Omega^*(d[d_1,d_2])}E^*(N;a,d[d_1,d_2]).\leqno(5.6)
$$
\medskip
Corresponding to $(4.8)$, we have
$$
\sum_d{\varrho(d)\over\varphi(d)}|\Omega^*(d)|
={W(R_0)\over V(R_0)}\left(1+O(e^{-\tau})\right),\quad V(R_0)
=\prod_{p\le R_0}\left(1-{1\over p}\right),\leqno(5.7)
$$
via the same reasoning. Also we have 
$$
\ca{T}^*=\sum_{|D|\le R}\Phi^*(D)\left(\sum_{D|K}\Delta^*(K)
\tilde\lambda_R(K;\ell)\right)^2,\leqno(5.8)
$$
where
$$
\Delta^*(D)=\prod_{P|D}
\left(\sum_{p\in P}{|\Omega^*(p)|\over p-1}\right),\leqno(5.9)
$$
and
$$
\Phi^*(D)={1\over\Delta^*(D)^2}
\prod_{P|D}\left(\sum_{p\in P}{|\Omega^*(p)|\over p-1}
\left(1-{|\Omega^*(p)|\over p-1}\right)\right).\leqno(5.10)
$$
Here we have actually $|\Omega^*(p)|=k-1$.
\medskip
We are about to show an effective lower bound of $\ca{T}^*$. We first note the
trivial inequality
$$
\ca{T}^*\ge \ca{T}^{**},\leqno(5.11)
$$
where the right side is the restriction of that of $(5.8)$ to $R/w\le |D|\le R$.
Inserting $(4.3)$ into $(5.8)$, we get
$$
\leqalignno{
\ca{T}^{**}=&\left({{\eufm S}(\ca{H})\over\ell!W(R_0)}\right)^2
\sum_{R/w\le|D|\le
R}\mu(D)^2\left({\Delta^*(D)\over\Delta(D)}\right)^2
{\Phi^*(D)\over\Phi(D)^2}&(5.12)\cr
\times&\left(\sum_{\scr{|K|\le R/|D|}\atop\scr{(K,D)=1}}
{\mu^2(K)\over\Phi(K)}\prod_{P|K}\left(1-{\Delta^*(P)\over\Delta(P)}
\right)\left(\log{R/|D|\over|K|}\right)^\ell\right)^2.
}
$$
This sum over $K$ can be handled with a simple modification
of the argument leading to $(3.17)$ besides employing $(3.19)$ with
an obvious change. In fact, we
may drop the condition
$K\in\ca{Y}(w)$, since $R/|D|\le w$. We have, 
for $\log R_0R_1\le \log y\ll\log R$,
$$
\sum_{\scr{|K|\le y}\atop\scr{P|K\Rightarrow P\nmid D,\, |P|\le y}}
{\mu^2(K)\over\Phi(K)}\prod_{P|K}\left(1-{\Delta^*(P)\over\Delta(P)}
\right)=V(R_0)\log y
\left(1+O\left({\log R_0\over\log y}
\right)\right), \leqno(5.13)
$$
uniformly in $D$.
The sum in question is then computed by integration by parts, and the result is
inserted into $(5.12)$ to give that
$$
\leqalignno{
&\qquad\ca{T}^{**}=\left({{\eufm
S}(\ca{H})\over(\ell+1)!}\cdot
{V(R_0)\over W(R_0)}\right)^2\left(1+O((\log\log
R)^{-1/5}))\right)&(5.14)\cr
\times&\left\{\sum_{\scr{|D|\le R}\atop
\scr{P|D\Rightarrow |P|\le w}}-
\sum_{\scr{|D|\le R/w}\atop
\scr{P|D\Rightarrow |P|\le w}}\right\}\mu(D)^2
\left({\Delta^*(D)\over\Delta(D)}\right)^2
{\Phi^*(D)\over\Phi(D)^2}\left(\log {R\over|D|}\right)^{2(\ell+1)}
}
$$
To estimate the part over $|D|\le R$, we proceed exactly as in $(4.12)$--$(4.15)$;
and the part over $|D|\le R/w$ as in $(4.13)$ or rather $(4.15)$,
appealing to $(3.17)$ and $(3.19)$ with an obvious change. In this way we
find that
$$
\ca{T}^{**}={{\eufm
S}(\ca{H})\over(k+2\ell+1)!}{2(\ell+1)\choose\ell+1}{V(R_0)\over W(R_0)}(\log
R)^{k+2\ell+1}\left(1+O(e^{-k\omega/3})\right),\leqno(5.15)
$$
which ends our treatment of the main term of $(5.4)$.
\bigskip
\noindent
{\bf 6.} We still need to consider the structure of $\ca{E}$, and it is
embodied in the assertion $(6.2)$ below.
\medskip
\noindent
{\bf Lemma 4.} {\it Under $(1.1)$, $(1.2)$, $(1.17)$, $(2.1)$, $(4.1)$, $(4.2)$,
$(4.3)$ $(4.6)$, $(4.7)$, $(4.9)$, it holds for any $h\in\ca{H}$ that
$$
\leqalignno{
&\sum_{N<n\le 2N}\varpi(n+h)\gamma(n;\ca{H})\tilde\Lambda_R(n;\ca{H},\ell)^2
&(6.1)\cr
&\ge{{\eufm S}(\ca{H})\over(k+2\ell+1)!}
{2(\ell+1)\choose\ell+1}N(\log R)^{k+2\ell+1}\left(1+O(e^{-k\omega/3})\right)
-\ca{E}_h(N;\ca{H}),
}
$$
as $N\to\infty$.
Here we have, for any $A, B\ge1$ such that $AB=R_0^{2\tau}R^{2+\omega}$, 
$$
\ca{E}_h(N;\ca{H})\le(\log N)^{2(k+\ell)+1}\sup_{\alpha,\,\beta}
\left|\sum_{a\le A,\,b\le B}\alpha_{a}\beta_{b}
\sum_{r\in\Omega^*(ab)}E^*(N;r,ab)\right|,\leqno(6.2)
$$
with $\Omega^*(p)=\Omega(p)\setminus\{-h\bmod p\}$ and $\Omega^*(p^v)=\emptyset$
$(v\ge2)$, where $\alpha,\,\beta$ run over vectors
such that  $|\alpha_{a}|\le1$, $|\beta_{b}|\le1$.
}
\medskip
\noindent
{\it Remark\/}. The above convention on $\Omega^*(ab)$ for non-square-free $ab$ can
in fact be replaced appropriately in practice. This is due to the fact that
in our construction below the situation $p^2|ab$ is possible only with $p\ge R_0$,
and the elimination of the contribution of those moduli is immediate. It should
also be stressed that we have in fact $\alpha_a=0$ or $1$, and $\beta_b=0$ or
$\varrho(d)$, $d\Vert b$, with $d$ being as in $(4.7)$.
\medskip
\noindent
The first term on the right of $(6.1)$ follows from $(5.4)$, $(5.7)$, $(5.11)$,
and $(5.15)$. As to $(6.2)$ we argue as follows: Returning to $(5.6)$, we
consider a generic pair $D_1, D_2$. Let $F$ be an arbitrary divisor of
$(D_1,D_2)$, the greatest common divisor of the pair. We restrict ourselves
to the situation in $(5.6)$ where $d_1\in D_1,\,d_2\in D_2$ and $(d_1,d_2)\in
F$. Let $D_1D_2/F=P_1P_2\cdots P_s$ with $|P_j|\le |P_{j+1}|$. Note that there
can be some $j$ such that $P_j=P_{j+1}$; in fact this is the case where $P_j$
divides $(D_1,D_2)/F$. We define 
$u$ to be such that $|P_1|\cdots|P_u|\le A$ but $|P_1|\cdots|P_{u+1}|> A$. It is
possible that there does not exist such $u$; then we are done. Otherwise, let 
$a\in P_1\cdots P_u$ and $a'\in P_{u+1}\cdots P_s$. 
Obviously we have $aa'\le R^2$. On the
other hand, we have $a\ge |P_1|\cdots|P_u|R_1^{-u}>A|P_{u+1}|^{-1}R_1^{-u}$,
because of the definition of the intervals given after $(2.1)$. Thus
$a'<R^{2+\omega}R_1^u/A$, as $|P_{u+1}|\le R^\omega$. Let $d$ be as in $(5.6)$,
and put $b=a'd$ we have $b<R^{2+\omega}R_0^\tau R_1^u/A<B$, since 
$u\ll (\log R)/\log R_0\ll(\log\log R)^{1/5}$. We are about to 
designate these $a,b$ as to be the same as in $(6.2)$; note that
$d[d_1,d_2]$ in $(5.6)$ are among the set of $ab$. Then we need to exclude those
$ab$ which are not square-free, for only those moduli are superfluous.
One way to employ here is to introduce a convention about $\Omega^*(p^v)$
$(v\ge2)$ as is done above. Finally, on noting $(4.5)$ as well as that the number
of triples $D_1,D_2, F$ is not larger than $\exp((\log\log R)^{9/10}\log 3)$, we
end the proof of
$(6.2)$.
\medskip
In the possible application to the problem about the gaps between primes,
we may assume that $k$ is large, and $\omega$ can be so small as 
$3(\log k)/k$. Hence the  size of $AB$ is essentially $R^{2+\varepsilon}$
with an arbitrarily small constant $\varepsilon>0$. With this, we see that
a combination of Lemmas 3 and 4 implies that 
if there exists a $C_1\ge 2(k+\ell+1)$ such that uniformly for $h$ in any
admissible
$\ca{H}$
$$
\ca{E}_h(N;\ca{H})\ll {N\over (\log N)^{C_1}},\quad \hbox{$R=N^{\theta/2}$
with an absolute constant $\theta>{1\over2}$},\leqno(6.3)
$$ 
then $(1.16)$ should follow. This hypothesis is
certainly less stringent than $(1.10)$ with $\theta>{1\over2}$. What is
interesting is that $(6.3)$ is true if the condition $r\in\Omega^*(ab)$ is
replaced  by $r\equiv r_0\pmod {ab}$ with a fixed integer $r_0$, as is proved in
[1]. It is, however, unclear how to extend the argument of [1] to the
situation with many residue classes as we require. 
\medskip
\noindent
{\it Concluding Remark\/}. The argument of our paper can be employed in a more
general setting:  With a {\it large} two-sided sifting density $\kappa$ (see, e.g.,
[5, p.\ 29]), the remainder term in the Selberg sieve admits a 
flexible bilinear form similar to the one proved by H. Iwaniec for Rosser's linear
sieve, although the level condition $MN\le D$, in the now common notation, has to
be replaced by the slightly weaker $MN\le D^{1+\delta}$ with $\delta\ll
(\log\kappa)/\kappa$, which is to be compared with $(6.2)$. In fact, this assertion
was obtained by the first author in early 1980's; however, any possible
application of it was not in his view then and even later when the relevant article
[7] was written.  He realised recently that his old method could be applied
to smoothing both Lemmas 1 and 2, and reached an earlier version of Lemmas 3 and 4.
Simultaneously and independently, the second author obtained the same.
\bigskip
\noindent
{\it Acknowledgements\/}. The present work is an outcome of the workshop
`{\it Gaps between Primes\/}' (November 28--December 3, 2005) at the American
Institute of Mathematics. The authors are indebted to the institute for the
financial aid that made it possible for them to attend the workshop.

\bigskip
\noindent
{\bf References}
\medskip
\noindent
\item{[1]} E. Bombieri, J.B. Friedlander and H. Iwaniec. 
Primes in arithmetic progressions to large moduli. Acta Math., {\bf 156}
(1986), 203--251.
\item{[2]} G.A. Goldston, S.W. Graham, J. Pintz, and C.Y. Y{\i}ld{\i}r{\i}m.
Small gaps between primes or almost 
primes. Preprint. arXiv:
math.NT/506067.
\item{[3]} D.A. Goldston, J. Pintz, and C.Y. Y{\i}ld{\i}r{\i}m. Small
gaps between primes II (Preliminary).  Preprint, February 8, 2005.
\item{[4]} G.A. Goldston, Y. Motohashi, J. Pintz, and C.Y. Y{\i}ld{\i}r{\i}m.
Small gaps between primes exist. Submitted.
\item{[5]} G. Greaves. Sieves in Number Theory. Ergebnisse Math., {\bf 43},
Springer Verlag, Berlin etc.\ 2001.
\item{[6]} Y. Motohashi. Sieve Methods and Prime Number Theory.
Tata IFR LN, vol.\ {\bf 72}, Springer Verlag, Berlin etc.\ 1983.
\item{[7]} Y. Motohashi. On the remainder term in the Selberg sieve. 
Number Theory in Progress - A. Schinzel Festschrift, Walter de Gruyter, 
Berlin $\cdot$ New York 1999, pp.\ 1053--1064.
\item{[8]} Y. Motohashi. An overview of the sieve method and its history.
Preprint (to appear in Sugaku Expositions, AMS). arXiv: math.NT/0505521. 

\bigskip
\noindent
{\srm
Yoichi Motohashi
\par
\noindent
Department of Mathematics,
Nihon University, 
\par
\noindent
 Surugadai, Tokyo 101-8308, JAPAN.
\par
\noindent
E-mail: ymoto@math.cst.nihon-u.ac.jp
\bigskip
\noindent
J\'anos Pintz
\par
\noindent
R\'enyi Mathematical 
Institute of the Hungarian Academy of Sciences,
\par\noindent
H-1364 Budapest, P.O.B. 127, HUNGARY
\par\noindent
E-mail: pintz@renyi.hu
}
\vfill\break
\vglue 1cm
\centerline{\title Addendum}
\vskip 0.7cm
We are indebted to Professor Terrence 
Tao for kindly pointing out
an inadequacy in the explanation 
following the formula $(5.14)$ of this paper of ours, 
designated as MP. The aim of the present addendum is
to provide a rectification to our argument. This
entails a replacement of the condition $(4.2)$ of MP on the
basic parameters $k,\ell,\omega$ by
a more restrictive one; see
$(**)$ below. It should be stressed, however, 
that being minor our correction 
does not make any difference to the main 
assertion of MP that the never-ending
appearances of bounded gaps between primes could
be established by assuming the possibility of
going past level $1\over 2$ in a mean prime number theorem
of the Bombieri--Vinogradov type in which the moduli are
restricted to smooth numbers, i.e., those $q$ whose prime
factors are all less than $q^\varpi$ with an arbitrary small
but fixed $\varpi>0$. We should repeat here the statement in
Lemma 4 of MP that our smoothing of the GPY sieve
induces a flexibility in the formulation of this conjectural 
mean prime number theorem, which allows for 
an appeal to the method due to Bombieri, Fouvry,
Friedlander and Iwaniec
concerning the distribution of primes in arithmetic 
progressions to large moduli. 
\par
As is widely known, Yitang Zhang
proved recently this conjectural part of MP, which
may now be called Zhang's mean prime number theorem; 
his argument depends on the works of the
four people above as we envisaged.
Therefore the bounded gaps between primes do occur 
infinitely often. It should be remarked that in Zhang's work
a smoothing of the GPY sieve is given as an essential
ingredient of his argument; but, including
its implication for the structure mentioned above of 
the mean prime number theorem to have been established,
his relevant reasoning is largely identical 
to the one we developed in MP
several years earlier. It appears to us that without
the smoothing of the GPY sieve 
the assertion on bounded gaps between primes
would not be deducible from Zhang's mean prime number
theorem.
\medskip
Now, turning to MP, immediately after $(5.14)$ it is
stated that (i) to estimate the part over $|D|\le R$, 
we proceed exactly as in $(4.12)$--$(4.15)$ and
(ii) the part over $|D|\le R/w$ as in $(4.13)$ 
or rather $(4.15)$,
appealing to $(3.17)$ and $(3.19)$ 
with an obvious change. 
\par
Neither (i) nor (ii) are quite adequate.
In order to rectify these, we note first the simple 
fact that $\Delta^*(P)\Phi^*(P)
> \Delta(P)\Phi(P)$ for any $P$, because of $(2.1)$. 
Namely, we may consider, instead of $(5.14)$, the expression
$$
\eqalign{
{\cal T}^{**}&\ge\left({{\eufm
S}(\ca{H})\over(\ell+1)!}\cdot
{V(R_0)\over W(R_0)}\right)^2\left(1+O((\log\log
R)^{-1/5}))\right)\cr
&\times\left\{\sum_{\scr{|D|\le R}\atop
\scr{P|D\Rightarrow |P|\le w}}-
\sum_{\scr{|D|\le R/w}\atop
\scr{P|D\Rightarrow |P|\le w}}\right\}{\mu(D)^2
\over\Phi^*(D)}
\left(\log {R\over|D|}\right)^{2(\ell+1)}.
}
$$
With this, the assertion
(i) is quite correct. As to (ii), we replace the factor 
$(\log R/|D|)^{2(\ell+1)}$ by
$$
\sum_{f=0}^{2(\ell+1)}{2(\ell+1)\choose f}(\log w)^f
\left(\log {R/w\over|D|}\right)^{2(\ell+1)-f},
$$
and consider the sum
$$
\sum_{|D|\le R/w}{\mu(D)^2\over\Phi^*(D)}
\left(\log {R/w\over|D|}\right)^{2(\ell+1)-f}.
$$
Note that the condition $P|D\Rightarrow |P|\le w$
has been dropped. This is equal to
$$
{(2(\ell+1)-f)!\over (k+2\ell+1-f)!}{W(R_0)\over
{\eufm S}(\ca{H})V(R_0)}\big(\log R/w\big)^{k+2\ell+1-f}
\left(1+O((\log\log R)^{-1/5}))\right).
$$
Hence we have
$$
\eqalign{
&\sum_{\scr{|D|\le R/w}\atop
\scr{P|D\Rightarrow |P|\le w}}{\mu(D)^2\over\Phi^*(D)}
\left(\log {R\over|D|}\right)^{2(\ell+1)}\cr
\le &\, Y(k,\ell;\omega)\cdot
{(2(\ell+1))!\over (k+2\ell+1)!}
{W(R_0)\over{\eufm S}(\ca{H})V(R_0)}
(\log R)^{k+2\ell+1}
\left(1+O((\log\log R)^{-1/5}))\right),
}
$$
with
$$
Y(k,\ell;\omega)=\sum_{f=0}^{2(\ell+1)}{k+2\ell+1\choose f}
\omega^f(1-\omega)^{k+2\ell+1-f};\quad
w=R^\omega.
$$
We have thus
$$
\ca{T}^{**}\ge{{\eufm
S}(\ca{H})\over(k+2\ell+1)!}{2(\ell+1)\choose\ell+1}
{V(R_0)\over W(R_0)}(\log R)^{k+2\ell+1}
\left(1-Y(k,\ell;\omega)+O(e^{-k\omega/3})\right).
$$
Denoting the $f$th term of $Y(k,\ell;\omega)$ 
by $A(f)$, we have
$$
A(f+1)/A(f)={k+2\ell+1-f\over f+1}{\omega\over 1-\omega}.
$$
This is decreasing as $f$ increases; and assuming
$$
{k\over 2(\ell+1)}{\omega\over 1-\omega}\ge1,\leqno(*)
$$
we have
$$
Y(k,\ell;\omega)\le(2\ell+3){k+2\ell+1\choose2(\ell+1)}
\omega^{2(\ell+1)}(1-\omega)^{k-1}.
$$
\par
Then, we make a drastic but practical specialisation:
We put, with an arbitrary constant $\alpha>0$,
$$
\omega=8{\alpha\over\sqrt{k}},\quad \ell=\alpha\sqrt{k},
\leqno(**)
$$
which amply satisfies $(*)$, provided $k$ is sufficiently
large. We find, by means of Stirling's formula, that
$$
Y(k,\ell;\varpi)<e^{-3k\omega/8}.
$$
Hence, we have obtained $(6.1)$ of MP
on $(**)$. 
\par
Ending our discussion, we 
remark that $(**)$ induces an obvious
change in the paragraphs of MP
following the proof of $(6.2)$;
however, this does not affect our overall conclusion
in the paper.
\bye